\documentclass{commat}

\usepackage{graphicx}

\newcommand{\HH}{\discretionary{-}{-}{-}}
\DeclareMathOperator{\diag}{diag}

\title[Associated Bessel functions related to one group of shifts]{A generalization of certain associated Bessel functions in connection with a group of shifts}

\author[J. Choi, I.A. Shilin]{Junesang Choi, Ilya A. Shilin}

\affiliation{
\address{J. Choi -- Department of Mathematics, Dongguk University, Gyeongju 38066, Republic of Korea}
\email{junesang@dongguk.ac.kr}
\address{I.A. Shilin -- Department of Higher Mathematics, National Research University MPEI,
Krasnokazarmennaya 14, Moscow 111250, Russia}
 \email{shilinia@mpei.ru}
\address{I.A. Shilin -- Department of Algebra, Moscow State Pedagogical University, Malaya Pirogovskaya 1, Moscow 119991, Russia}
 \email{ilyashilin@li.ru}
}



\abstract{Considering the kernel of an integral operator
intertwining two realizations of the group of motions of the
pseudo-Euclidean space, we derive two formulas for series
containing Whittaker's functions or Weber's parabolic cylinder
functions. We can consider this kernel as a special function.
Some particular values of parameters involved in this special function are found to coincide with
certain variants of Bessel functions. Using these connections, we
also establish some analogues of orthogonality relations for
Macdonald and Hankel functions.
}

\firstpage{103}

\VOLUME{30}

\DOI{https://doi.org/10.46298/cm.9305}

\begin{paper}

\section{Introduction and preliminaries}\label{sec1}

It is well known that any Lie group depends on a finite set of
continuously changing parameters. The cardinality of this set is
small for a group of a low dimension. In~this case, the matrix
elements of representation operators, the matrix elements of~bases
transformations, kernels of the corresponding integral operators, intertwined different realizations of
representations, can be expressed in terms of classical special
functions. For more complicated groups, the above matrix elements
and kernels of subrepresentations to some subgroups and kernels of
intertwining operators are found to yield new special functions, which can
be considered either generalizations or analogues of known
(classical) special functions. For instance, in~\cite{sc}, the matrix elements of the restriction of the
representation of Lorentz group on some diagonal matrices were shown to be able to be
expressed in terms of modified hyper Bessel functions of the first
kind. The connection between these matrix elements written in two
different bases of a representation space leads to new formulas
for series containing above-mentioned hyper functions and converging to
(ordinary) modified Bessel functions.

Vilenkin~\cite{v} showed that many known and new
properties of variants of Bessel functions are related to representations of
the group $ISO(1,1)$ (denoted by $MH(2)$ there) of motions
of the pseudo-Euclidean plane and group $M(n)$ of motions of
Euclidean $n$-dimensional space. In this paper, we consider an
one-parameter subgroup in more complicated group $ISO(n,1)$ and show
that the kernel of an intertwined integral operator can be
considered as a generalization of Macdonald and Hankel functions, in a sense that
some simple cases of the kernel coincide with those variants of Bessel functions.

Recall that the pseudo-Euclidean ($n+1$)-dimensional space of
signature
\[
\{-,\dots,-,+\}
\]
is the linear space in
$\mathbb{R}^{n+1}$ endowed with the bilinear form
\[
\mu(x,y) = x_{n+1}y_{n+1}-\sum\limits_{i = 1}^n x_iy_i\,.
\]
Here and
throughout, let $\mathbb{N}$, $\mathbb{Z}$, $\mathbb{R}$,
$\mathbb{R}^+$, and $\mathbb{C}$ denote the sets of positive
integers, integers, real numbers, positive real numbers, and
complex numbers, respectively, and let $\mathbb{N}_0: = \mathbb{N}
\cup \{0\}$ and $\mathbb{R}^+_0: = \mathbb{R}^+ \cup \{0\}$. Let
$\tilde\mu(x) = \mu(x,x)$ be the corresponding quadratic form. A
motion of the pseudo-Euclidean space is an isometry with respect
to the distance $\sqrt{\tilde\mu(x)}$ that preserves orientation.
These motions form a group denoted by $ISO(n,1)$. For any $g\in
ISO(n,1)$ and any $x\in\mathbb{R}^{n+1}$, we have
\[
g(x) = g_0x+g_1\,,
\]
where $g_0\in SO(n,1)$ and $g_1$ is a shift vector. Therefore,
$g = g(g_0,g_1)$. Let $\mathfrak{D}$ be the linear space consisting
of infinitely differentiable functions defined on the upper cone
$\tilde\mu(x) = 1$, where $x = \left(x_1,\ldots,x_{n+1}\right)$ with
each $x_j \in \mathbb{R}^+$.

For any $\tau \in \mathbb{C}$, we consider a map $T_\tau$ defined
by
\begin{equation*}
  T_\tau\colon ISO(n,1)\rightarrow
GL(\mathfrak{D}), \quad
f\mapsto\exp\big(-\tau\mu(g_1,x)\big)\,f(g_0^{-1}x)
\end{equation*}
where $GL(\mathfrak{D})$ is the multiplicative group of linear
operators of $\mathfrak{D}$ whose rank is nonzero. Since
for any $g,\tilde g\in ISO(n,1)$ and any $x\in\mathbb{R}^{n+1}$ we
have
\[
\tilde gg(x) = \tilde g(g_0x+g_1) = \tilde g_0g_0x+\tilde
g_0g_1+\tilde g_1\,.
\]
So we obtain $\tilde gg = \hat g(\hat g_0,\hat
g_1)$, where $\hat g_0 = \tilde g_0g_0$ and $\hat g_1 = \tilde
g_0g_1+\tilde g_1$.
Then
\begin{align*}[T_\tau(\tilde g)T_\tau(g)]
(f(x)] &= [T_\tau(\tilde g)]\bigl((\exp\big(-\tau\mu(g_1,x)\big)\,f(g_0^{-1}x)\bigr)\\
 &= \exp\bigl(-\tau,\mu(\hat
g_1,x)\bigr)\,f(\hat g_0^{-1}x)\,.
\end{align*}
We thus find that $T_\tau$ is a homomorphism.

In order to simplify the representation of the group $SH(n+1)$ of
hyperbolic rotations in $\mathbb{R}^{n+1}$, Vilenkin~\cite[Chapter
10]{v} employed the so-called horosphere method. Indeed, Vilenkin
used the Gelfand-Graev integral transformation, which maps
$f\in\mathfrak{D}$ into the space $L$ of $\hat\sigma$-homogeneous
functions defined on the intersection of the cone $\tilde\mu(x) = 0$
and the plane $x_{n+1} = 1$. Let
$K = (k_0,\ldots,k_{n-2})\in\mathbb{Z}^{n-1}$, $k_0\geqslant
k_1\ldots\geqslant k_{n-3}\geqslant|k_{n-2}|$, $C_k^r$ be the
Gegenbauer polynomials, and
\begin{align*}
{(A_K)}^2 &= 2^{n-3}\pi^{-\frac12}\Gamma^{-1}\Bigl(\frac{n}2\Bigr)\\
&\quad \times
\prod\limits_{i = 0}^{n-3}\frac{2^{2k_{i+1}-i}(k_i-k_{i+1})!(2k_i+n-2-i)
\Gamma^2 \bigl(\frac{n-i}2+k_{i+1}-1\bigr)}{\Gamma(k_i+k_{i+1}+n-i-2)}\,.
\end{align*}
It was shown that for any $f\in\mathfrak{D}$, the `coordinates'
$a_K$ of the image of $f$ (arising after action of Gelfand-Graev
transform) with respect to the basis
\[
\Xi_K(x) = A_K\, \prod\limits_{i = 0}^{n-3}C_{ki-k_{i+1}+1}^{\frac{n-i-2}2}(\cos\varphi_{n-i-1})\sin^{k_i+1}
\varphi_{n-i-1}\exp(\pm\mathbf{i}k_{n-2}\varphi_1)
\]
(where $A_k$ is the normalizing factor) can be expressed as the integral transform (see~\cite[Entry 10.5.4]{v})
\[
a_K = \mathsf{I}[f](\sigma) = \int\limits_{\tilde\mu(x) = 1}\ker(x,K,\sigma)f(x)\,\mathrm{d}x\,.
\]
Here the kernel $\ker(x,K,\sigma)$ coincides with the matrix elements of the
$SH(n+\nobreak1)$\HH representation in the space $L$ situating in the `zero'
column which is exactly described by formula~\cite[Entry
10.4.9]{v}. The corresponding formulas for the inverse transform
are given by (see~\cite[Entries 10.5.5 and 10.5.6]{v})
\begin{equation}\label{imp}
f(x) = \frac{{(-1)}^{\frac{n-\mathbf{\eta}}2}\,\mathbf{i}}{2^{n}\pi^{\frac{n}2}\Gamma\bigl(\frac{n}2\bigr)}\sum\limits_{K}\int\limits_{c-{\bf
i}\infty}^{c+\mathbf{i}\infty}(\sigma)_{n-1}\,[\cot(\pi\sigma)]^{1-\mathbf{\eta}}\,a_K(\sigma)\,\ker (x,K,1-n-\sigma)\,\mathrm{d}\sigma\,.
\end{equation}
Here and elsewhere, $\mathbf{i} = \sqrt{-1}$
and $\mathbf{\eta}$ is the remainder of $n$ divided by $2$, that
is, either $\mathbf{\eta} = 0$ or $\mathbf{\eta} = 1$. We also recall some
functions and notations, which are used in the following sections.
The generalized hypergeometric series $_pF_q$ $\left(p,\,q \in
\mathbb{N}_0\right)$ is defined by (see~\cite[p.~73]{Rain}):
\begin{align*}\label{pFq}
 {}_pF_q \left[ \aligned \alpha_1,\dots,\alpha_p &;\\
                       \beta_1,\dots,\beta_q &; \endaligned
                   z \right] &= \sum_{n = 0}^\infty \frac{(\alpha_1)_n \dots (\alpha_p)_n
}{(\beta_1)_n \dots (\beta_q)_n} \frac{z^n}{n!} \\
               &= {}_pF_q (\alpha_1,\dots,\alpha_p; \beta_1,\dots,\beta_q;z)\,,
\end{align*}
where $(\lambda)_\nu$ denotes the Pochhammer symbol which is defined (for $\lambda,\nu \in \mathbb{C}$), in~terms of the familiar Gamma function $\Gamma$, by
\begin{equation*}\label{Poch-symbol}
  (\lambda)_\nu : = \frac{\Gamma (\lambda +\nu)}{\Gamma (\lambda)}
 =
\begin{cases}
1 & (\nu = 0; \lambda \in \mathbb{C}\setminus \{0\}) \\
         \lambda (\lambda +1) \dots (\lambda+n-1) & (\nu = n \in \mathbb{N};\lambda \in \mathbb{C}),
   \end{cases}
\end{equation*}
it being understood conventionally that $(0)_0: = 1$. The Gegenbauer function $C_\nu^\lambda(z)$ is defined by (see, e.g.,~\cite[p. 791]{v3})
  \begin{equation*}\label{GegenFt}
    C_\nu^\lambda(z) = \frac{\Gamma (2\lambda +\nu)}{\Gamma (2\lambda)\,\Gamma (\nu+1)}\,
      {}_2F_1 \Bigl(-\nu,2\lambda +\nu; \lambda +\frac{1}{2};\frac{1-z}{2}\Bigr)\,,
  \end{equation*}
the cases $\nu = n \in \mathbb{N}_0$ of which are the Gegenbauer polynomials. The associated Legendre function of the first kind $P_\nu^\mu(z)$
is defined by (see, e.g.,~\cite[p. 795]{v3})
\begin{gather}\label{ALF-1st}
  P_\nu^\mu(z) = \frac{1}{\Gamma (1-\mu)}\Bigl{(\frac{z+1}{z-1}\Bigr)}^{\frac{\mu}{2}}
    {}_2F_1 \Bigl(-\nu,\,\nu+1;1-\mu;\frac{1-z}{2}\Bigr) \\
(\lvert\arg (z\pm 1)\rvert<\pi;\quad \mu \in \mathbb{C} \setminus \mathbb{N})\,. \notag
\end{gather}
 The confluent hypergeometric function $\Psi$ is given by (see, e.g.,~\cite[Entry 13.6.21]{Olver2010-NIST})
\begin{equation*}\label{U-2F0}
  \Psi(a,b; z) = z^{-a} {}_2F_0 \Bigl(a, a-b+1;-; - \frac{1}{z} \Bigr)\,.
\end{equation*}
The Whittaker function of the second kind $W_{\kappa,\mu}(z)$ is defined by (see, for instance, \cite[p.~264, Eq.~(5)]{EMOT-I}
and~\cite[p. 334, Entry 13.14.3]{Olver2010-NIST})
\begin{equation}\label{WhF-2nd}
 W_{\kappa,\mu}(z) = e^{-\frac{z}{2}}\, z^\kappa\,{}_2F_0 \Bigl(\frac{1}{2}-\kappa+\mu,\frac{1}{2}-\kappa-\mu;-;-\frac{1}{z}\Bigr)
\end{equation}
\begin{equation*}
  \Bigl(\lvert\arg (z)\rvert<\pi, |z|>0;\quad \frac{1}{2}-\kappa+\mu, \frac{1}{2}-\kappa-\mu \in \mathbb{C} \setminus \mathbb{Z}_0^- \Bigr)\,.
\end{equation*}
Recall the following relation (see, e.g.,~\cite[Entry 13.14.3]{Olver2010-NIST})
\begin{equation}\label{U-Whit-2nd}
  W_{\kappa,\mu}(z) = z^{\mu+\frac12}\,e^{-\frac{z}2}\,\Psi\Bigl(
\frac12-\kappa+\mu,\,2\mu+1;z\Bigr)\,.
\end{equation}
The parabolic cylinder function $D_\nu (z)$ is given by (see, for instance, \cite[p.~674]{Bryc} and~\cite[p.~792]{v3})
\begin{equation*}\label{PCF}
  D_\nu (z) = 2^{\frac{\nu}{2}} e^{-\frac{z^2}{4}} \Psi \Bigl(-\frac{\nu}{2}, \frac{1}{2}; \frac{z^2}{4}\Bigr)\,.
\end{equation*}
The modified Bessel function of the $1$st kind $I_\nu(z)$ is given by (see, e.g.,~\cite[p.~794]{v3})
\begin{equation*}\label{MBF-1st}
  I_\nu(z) = \frac{1}{\Gamma (\nu+1)} \Bigl{(\frac{z}{2}\Bigr)}^\nu\, {}_0F_1 \Bigl(-;\nu+1; \frac{z^2}{4} \Bigr) = e^{- \frac{\nu \pi \textbf{i}}{2}} J_\nu (e^{\frac{\pi \textbf{i}}{2}}z)\,,
  \end{equation*}
where $J_\nu (z)$ is the Bessel function of the $1$st kind given by (see, e.g.,~\cite[p. 794]{v3})
\begin{equation*}\label{BF-1st}
  J_\nu (z) = \frac{1}{\Gamma (\nu+1)} \Bigl{(\frac{z}{2}\Bigr)}^\nu {}_0F_1 \Bigl(-;\nu+1;- \frac{z^2}{4} \Bigr)\,.
\end{equation*}
The MacDonald function (modified Bessel function of the $3$rd kind) $K_\nu(z)$ is given by (see, e.g.,~\cite[p. 794]{v3})
\begin{equation}\label{MacF}
 K_\nu(z) = \frac{\pi [I_{-\nu}(z) -I_{\nu}(z) ]}{2\, \sin (\nu \pi)}\,, \quad (v \in \mathbb{C} \setminus \mathbb{Z}),
\end{equation}
whose integral representation of pure imaginary index $\nu = \textbf{i} t$ and the argument $z = x \in \mathbb{R}$ was used as a definition in~\cite[p. 873, Eq. (16)]{cr}:
\begin{equation*}\label{MacF-purei}
  K_{\textbf{i} t}(x) = \int\limits_{0}^{\infty} e^{-x\cosh u} \cos (tu)\,\textrm{d}u\,.
\end{equation*}
The Hankel functions of the first and second kind (the Bessel functions of the third kind) are given, respectively, by
(see, e.g.,~\cite[p. 675]{Bryc})
\begin{equation*}\label{HFts}
  H_\nu^{(1)}(z) = J_\nu (z) + \textbf{i}\, Y_\nu (z) \quad \text{and} \quad H_\nu^{(2)}(z) = J_\nu (z) - \textbf{i}\, Y_\nu (z),
\end{equation*}
where $Y_\nu (z)$ is the Bessel function of the second kind (the Neumann function) (see, e.g.,~\cite[p. 677]{Bryc})
\begin{equation*}\label{BF-2nd}
 Y_\nu (z) = \frac{\cos (\nu \pi)\,J_\nu (z)-J_{-\nu} (z)}{\cos (\nu \pi)}\,, \quad (\nu \in \mathbb{C} \setminus \mathbb{Z}).
\end{equation*}

The Kronecker symbol $\delta_{m,n}$ is defined by $\delta_{m,n} = 1$ when $m = n$ and $\delta_{m,n} = 0$ when $m \ne n$.

\section{Kernels of integral operators of the subrepresentation to
the subgroup of shifts along the axis $Ox_{n+1}$}\label{sec2}

Consider the subgroup
\[
H = \{h(\lambda): = g(\diag(1,\ldots,1),(0,\ldots,0,\lambda))\mid\lambda\in\mathbb{R}\}
\]
in $M(n,1)$. In this section and elsewhere, we deal with the
integral operator $\mathsf{B}$ which acts in the `space of functions
$a_K(\sigma)$' and corresponds with the representation operator
$T_\tau\big(h(\lambda)\big)$, where $h(\lambda)\in H$.

\begin{lemma}\label{lem-BI}
For any $h(\lambda)\in H$, the integral operator $\mathsf{B}$ can be
written as a Barnes integral
\[
\mathsf{B}[a_K(\sigma)] = \int\limits_{c-\mathbf{i}\infty}^{c+\mathbf{i}\infty}\widetilde{\ker}(\tau, \lambda;
n,K,\sigma,\hat\sigma)a_{\hat K}(\hat\sigma)\,\mathrm{d}\hat\sigma\,,
\]
where the kernel admits the integral
representation
\begin{multline}\label{ir}
\widetilde{\ker}(\tau, \lambda;n,K,\sigma,\hat\sigma) = \bigl[{(-1)}^{\frac{n}2}\cot(\pi\hat\sigma)\delta_{\mathbf{\eta},0}+{(-1)}^{\frac{n+1}2}\delta_{\mathbf{\eta},1}\bigr]2^{-1}\mathbf{i}(-\sigma)_{k_0}\\
 \times (\hat\sigma)_{n-1}(\hat\sigma+n-1)_{k_0}
\int\limits_1^{\infty}P_{\frac{n}2+\sigma-1}^{1-k_0-\frac{n}2}(t)\,P_{-\frac{n}2-\hat\sigma}^{1-k_0-\frac{n}2}(t)
\,e^{-\tau\lambda t}\,\mathrm{d}t.
\end{multline}
\end{lemma}

\begin{proof}
The inverse transform formulas~\eqref{imp} (which will be applied in this proof) coincide
for even and odd cases up to factor $\pm\cot(\pi\sigma)$.
Therefore, it is sufficient to prove the result for an arbitrary
odd $n$. Introduce the following parametrization on
$\tilde\mu(x) = 1$, where $x = (x_1,\dots,x_{n+1})$ with
$x_{n+1} \in \mathbb{R}^+$ (see~\cite[Entry 10.1.1]{v}):
\begin{gather*}
x_i = \sum\limits_{j = 1}^n \delta_{i,j}\sinh\theta_n\,\prod\limits_{s = i}^{n-1}\sin\theta_s\sum\limits_{t = 2}^n \delta_{t,i}\cos\theta_{i-1}
+\delta_{i,n+1}\cosh\theta_n, \\
\theta_1\in[-\pi,\pi], \theta_2,\dots,\theta_{n-1}\in[0,\pi], \theta_n \in \mathbb{R}^+_{0}.
\end{gather*}
Given $x = (x_1,\dots,x_{n+1})$ with $x_{n+1} \in \mathbb{R}^+$, the corresponding $SH(n+1)$-invariant measure on $\tilde\mu(x) = 1$ is given by (see~\cite[Entry 10.1.6]{v})
 \[
\mathrm{d}x = \sinh^{n-1}\theta_n\prod\limits_{i = 1}^{n-1}\sin^{i-1}\theta_i\,\mathrm{d}\theta_i\,.
\]
Then we find from~\cite[Entry 10.4.9]{v} that
\begin{align*}
\ker{}&{}(x,K,\sigma)\\
={ }&{ } {(-1)}^{k_0}2^{\frac{n}2-1}\Gamma(\sigma+1)\Gamma^{-1}(\sigma-k_0+1)\sinh^{1-\frac{n}2}\theta_n \exp(\pm\mathbf{i}k_{n-2}\theta_1) \\
&{ } \times \Bigl[\Gamma\Bigl(\frac{n}2\Bigr)\prod\limits_{i = 0}^{n-3}\frac{2^{2k_{i+1}+n-i-4}(n+2k_i-i-2) (k_i-k_{i+1})!\Gamma\bigl(\frac{n-i}2+k_{i+1}-1\bigr)}{\pi^{\frac{n}2-1}\Gamma(n+k_i+k_{i+1}-i-2)} \Bigr]^{\frac12} \\
&{ } \times P_{\frac{n}2+\sigma-1}^{1+k_0-\frac{n}2}(\cosh\theta_n) \prod\limits_{i = 0}^{n-3}C_{k_i-k_{i+1}}^{\frac{n-i}2+k_{i+1}-1}(\cos\theta_{n-i-1}) \sin^{k_{i+1}}\theta_{n-i-1}
\end{align*}
and
$[T_\tau(h)][f(x)] = \exp(-\tau\lambda\cosh\theta_n)f(x)$.
Therefore we obtain
\begin{align*}
\mathsf{B}[a_k(\sigma)]
={ }&{ }\frac{{(-1)}^{\frac{n}2}\,\mathbf{i}2^{1-n}}{\pi^{\frac{n}2}\Gamma\bigl(\frac{n}2\bigr)}\sum\limits_{\hat
K}\int\limits_{c-\mathbf{i}\infty}^{c+\mathbf{i}\infty}(\hat\sigma)_{n-1}\,a_{\hat K}(\hat\sigma)\,\mathrm{d}\hat\sigma\\
&{ }\times \ker (x,K,\sigma)\ker (x,\hat
K,1-n-\hat\sigma)\int\limits_{-\pi}^\pi\exp\bigl[\pm\mathbf{i}(k_{n-2}-\hat k_{n-2})\theta_1\bigr]\,\mathrm{d}\theta_1\\
&{ } \times \prod\limits_{i = 0}^{n-3}\int\limits_0^\pi
C_{k_i-k_{i+1}}^{\frac{n-i}2+k_{i+1}-1}(\cos\theta_{n-i-1})
C_{\hat k_i-\hat k_{i+1}}^{\frac{n-i}2+\hat
k_{i+1}-1}(\cos\theta_{n-i-1})\\
&{\quad} \times \sin^{k_{i+1}+\hat
k_{i+1}+n-i-2}\theta_{n-i-1}\,\mathrm{d}\theta_{n-i-1}\\
&{ }\times \int\limits_0^{+\infty}P_{\frac{n}2+\sigma-1}^{1-k_0-\frac{n}2}(\cosh\theta_n) P_{-\frac{n}2-\hat\sigma-1}^{1-\hat
k_0-\frac{n}2}(\cosh\theta_n)\,\sinh\theta_n\,\mathrm{d}\theta_n.
\end{align*}
Obviously $\mathsf{B}[a_k(\sigma)] = 0$ for $k_{n-2}\ne\hat k_{n-2}$. Otherwise, in view of
orthogonality property for Gegenbauer polynomials (see, e.g.,~\cite[p. 198]{bell})
\[
\int\limits_{-1}^1 {(1-x^2)}^{\varrho-\frac12}\,C_k^\varrho(x)\,C_m^\varrho(x)\,\mathrm{d}x = \frac{2^{1-2\lambda}\pi\Gamma(k+2\varrho)\delta_{k,m}}{k!(k+\varrho)\Gamma^2 (\varrho)}\,,
\]
we see that the integral
\[
\int\limits_0^\pi
C_{k_{n-3}-k_{n-2}}^{\frac12+k_{n-2}}(\cos\theta_2) C_{\hat
k_{n-3}-\hat k_{n-2}}^{\frac12+\hat k_{n-2}}(\cos\theta_2)
\sin^{k_{n-2}+\hat k_{n-2}+1}\theta_2\,\mathrm{d}\theta_2 = 0
\]
 for $k_{n-3}\ne\hat k_{n-3}$. Since, for $K = \hat K$
\[
\int\limits_{-\pi}^\pi\exp\bigl[\pm\mathbf{i}(k_{n-2}-\hat k_{n-2})\theta_1\bigr]\,\mathrm{d}\theta_1 = 2\pi
\]
and
\begin{multline*}
\prod\limits_{i = 0}^{n-3}\int\limits_0^\pi
C_{k_i-k_{i+1}}^{\frac{n-i}2+k_{i+1}-1}(\cos\theta_{n-i-1})
C_{\hat k_i-\hat k_{i+1}}^{\frac{n-i}2+\hat
k_{i+1}-1}(\cos\theta_{n-i-1})\\
 \times \sin^{k_{i+1}+\hat
k_{i+1}+n-i-2}\theta_{n-i-1}\,\mathrm{d}\theta_{n-i-1} = \pi^{n-2}\Bigl[\Gamma\Bigl(\frac{n}2\Bigr){(A_K)}^2 \Bigr]^{-1}\,,
\end{multline*}
we have
\begin{align*}
\mathsf{B}[a_k(\sigma)] &= {(-1)}^{\frac{n+1}2}\,2\mathbf{i}\Gamma(\sigma+1)\Gamma^{-1}(1+\sigma-k_0)\\
&\quad \times \int\limits_{c-\mathbf{i}\infty}^{c+\mathbf{i}\infty}(\hat\sigma)_{n-1}{(2-n-\hat\sigma)}^{k_0}a_K(\hat\sigma)\,\mathrm{d}\hat\sigma\\
&\quad \times \int\limits_1^{+\infty}P_{\frac{n}2+\sigma-1}^{1-k_0-\frac{n}2}(t)
P_{-\frac{n}2-\hat\sigma}^{1-k_0-\frac{n}2}(t)\exp(-\tau\lambda
t)\,\mathrm{d}t\,.
\end{align*}
This completes the proof.
\end{proof}

We find from Lemma~\ref{lem-BI}
that the kernel $\widetilde{\ker}(\tau, \lambda;n,K,\sigma,\hat\sigma)$ does not depend on
$k_1,\dots,k_{n-2}$, in fact,
\[
\widetilde{\ker}(\tau,
\lambda;n,K,\sigma,\hat\sigma)\equiv\widetilde{\ker}(\tau,
\lambda;n,k_0,\sigma,\hat\sigma)\,.
\]

\begin{lemma}\label{lem-Ker}
 The kernel $\widetilde{\ker}(\tau, \lambda;n,k_0,\sigma,\hat\sigma)$ admits the following
series representation:
\begin{itemize}
\item The case $\mathbf{\eta} = 0$
 \begin{align}
&{}\widetilde{\ker}(\tau, \lambda;n,k_0,\sigma,\hat\sigma) \\
&{\,} = {(-1)}^{\frac{n}2} {(2\tau\lambda)}^{-\frac{n+\sigma+\hat\sigma}2}\mathbf{i}\cot(\pi\hat\sigma)\nonumber \\
&{\ }\times (-\sigma)_{k_0}(\hat\sigma+n-1)_{k_0}(\hat\sigma)_{n-1}\left[\Gamma\left(\frac{n}2+k_0\right)\right]^{-1}\nonumber \\
&{\ } \times \sum\limits_{i = 0}^\infty\frac{\left(1-\sigma-\frac{n}2\right)_i(k_0-\sigma)_i}{i!}{}_4F_3\left[
\begin{array}
{r} 1-\frac{n}2-k_0-i,1+\hat\sigma-\frac{n}2,k_0-\hat\sigma,-i;\\
\frac{n}2+k_0,\frac{n}2+\sigma-i,1+\sigma-k_0-i;
\end{array}
1\right]\nonumber \\
&{\ } \times W_{\frac{\sigma+\hat\sigma}2-k_0-i,\frac{n+\sigma+\hat\sigma-1}2}(2\tau\lambda). \label{sr1}
\end{align}

\item The other case $\mathbf{\eta} = 1$
\begin{align*}
  \widetilde{\ker}(\tau,
\lambda;n,k_0,\sigma,\hat\sigma) &= {(-1)}^{\frac{n+1}2}
{(2\tau\lambda)}^{\frac{\hat\sigma-\sigma-1}2}\mathbf{i} (-\sigma)_{k_0}(\hat\sigma+n-1)_{k_0}\\
&\quad \times (\hat\sigma)_{n-1}\left[\Gamma\left(\frac{n}2+k_0\right)\right]^{-1} \sum\limits_{i = 0}^\infty\frac{\left(1-\sigma-\frac{n}2\right)_i(k_0-\sigma)_i}{i!}\\
&\quad \times{}_4F_3\left[
\begin{array}
{r}
1-\frac{n}2-k_0-i,\frac{n}2+\hat\sigma,n+k_0+\hat\sigma-1,-i;\\
\frac{n}2+k_0,\frac{n}2-\sigma-i,1+\sigma-k_0-i;
\end{array}
1\right]\\
&\quad \times
W_{\frac{1+\sigma-\hat\sigma-n}2-k_0-i,\frac{\sigma-\hat\sigma}2}(2\tau\lambda)\,.
\end{align*}
\end{itemize}
\end{lemma}

\begin{proof}
We prove this theorem only for the case $\mathbf{\eta} = 1$.
Using~\eqref{ALF-1st} and a product formula of functions ${}_pF_q$ (see, e.g.,~\cite[p. 441, Entry 7.2.3-44]{v3})
\[
{}_pF_q\left[
\begin{array}
{r}a_1,\dots,a_p;\\
b_1,\dots,b_q;
\end{array}
cz\right]
{}_rF_s\left[
\begin{array}
{r}\hat a_1,\dots,\hat a_r;\\
\hat b_1,\dots,\hat b_s;
\end{array}
\hat cz\right] = \sum\limits_{i = 0}^\infty\gamma_i z^i \,,
\]
 where
\[
\resizebox{\textwidth}{!}{$
\gamma_i = \frac{c^i \prod\limits_{j = 1}^p (a_j)_i}{i!\prod\limits_{j = 1}^q (b_j)_i}{}_{q+r+1}F_{p+s}
\left[
\begin{array}
{r}
-i,1-b_1-i,\dots,1-b_q-i,\hat a_1,\dots,\hat a_r;\\
1-a_1-i,\dots,1-a_p-i,\hat b_1,\dots,\hat b_s;
\end{array}
\frac{{(-1)}^{p+q+1}\hat c}c
\right],
$}
\]
we obtain
\begin{align}
&\int\limits_1^{\infty}P_{\frac{n}2+\sigma-1}^{1-k_0-\frac{n}2}(t) P_{-\frac{n}2-\hat\sigma}^{1-k_0-\frac{n}2}(t)\exp(-\tau\lambda t)\,\mathrm{d}t \\
&{\ } = \frac{2^{1-\sigma-\hat\sigma}}{\Gamma^2 \bigl(\frac{n}2+k_0\bigr)}\nonumber\\
&{\ \ } \times \sum\limits_{i = 0}^\infty \frac{\left(1-\sigma-\frac{n}2\right)_i\,(k_0-\sigma)_i}{i!\,\left(\frac{n}2+k_0\right)_i}\, {}_4F_3\left[
\begin{array}
{l}-i,\hat\sigma+n+k_0-1,\hat\sigma+\frac{n}2,1-k_0-i-\frac{n}2\\
\frac{n}2+k_0,\frac{n}2+\sigma-i, 1+\sigma-k_0-i
\end{array}
\,1\right]\nonumber\\
&{\ \ } \times \exp(-\tau\lambda)\,\int\limits_0^{\infty}t^{\frac{n}2+k_0+i-1}\,{(t+2)}^{\sigma-\hat\sigma-k_0-i-\frac{n}2}\,\exp(-\tau\lambda t)\,\mathrm{d}t. \label{b}
\end{align}
Using a known integral formula (see, e.g.,~\cite[Entry 2.3.6.-9)]{v1})
\begin{equation*}\label{IF-1}
  \int\limits_0^{\infty}t^{\alpha-1}{(t+u)}^{-\beta}e^{-pt}\,\mathrm{d}t = \Gamma(\alpha)u^{\alpha-\beta}U(\alpha;1+\alpha-\beta;up)
\end{equation*}
\begin{equation*}
  \left(\lvert\arg(u)\rvert<\pi, \min\{\Re(\alpha),\Re(p)\}>0\right)
\end{equation*}
and the relation~\eqref{U-Whit-2nd},
we complete the proof for the case $\mathbf{\eta} = 1$. Similarly the case $\mathbf{\eta} = 0$
can be shown.
\end{proof}

We find from Lemma~\ref{lem-Ker}
 that the kernel $\widetilde{\ker}(\tau,
\lambda;n,k_0,\sigma,\hat\sigma)$ depends on the product
$u: = \tau\lambda$, in fact,
\[
\widetilde{\ker}(\tau,
\lambda;n,k_0,\sigma,\hat\sigma)\equiv\widetilde{\ker}(u;n,k_0,\sigma,\hat\sigma)\,.
\]

It is noted that some particular values of parameters of $\widetilde{\ker}(u;n,k_0,\sigma,\hat\sigma)$ can
yield certain other series representations. For example, replacing
$\hat\sigma$ by $\sigma+\frac12$ in~\eqref{b}, we can use a known formula (see, e.g.,~\cite[Entry 2.3.6.-12]{v1})
\begin{equation}\label{c}\int\limits_0^{\infty}t^{\alpha-1}{(x+u)}^{-\alpha-\frac12}e^{-pt}\,\mathrm{d}t = 2^\alpha u^{-\frac12}\Gamma(\alpha)e^{\frac{pu}2}D_{-2\alpha}\left(\sqrt{2pu}\right)
\end{equation}
\begin{equation*}
  (\lvert\arg(u)\rvert<\pi,\min\{\Re(\alpha),\Re(p)\}>0)\,.
\end{equation*}
Then we may obtain the series involving the product of ${}_4F_3$ and the
parabolic cylinder function.

\section{Series representations of Macdonald functions}\label{sec3}

 Using the result in Lemma~\ref{lem-Ker}, the Macdonald function in~\eqref{MacF}
 can be expressed as a series involving the Whittaker function of the second kind $W_{\kappa,\mu}(z)$~\eqref{WhF-2nd},
 asserted in~the following theorem.

\begin{theorem}\label{thm1-sec3}
Let $n$ be even and $\Re(u)>0$. Then
\begin{equation}\label{thm1-sec3-eq1}
 K_{\sigma+\frac{n-1}2}(u) = {(2u)}^{-\frac{\sigma+1}2}\pi^{\frac12}\sum\limits_{i = 0}^\infty
\frac{(-\sigma)_i\left(1-\sigma-\frac{n}2\right)_i}{i!}W_{\frac\sigma2-i,\frac{\sigma+n-1}2}(2u)\,.
\end{equation}
\end{theorem}

\begin{proof}
Considering that,
for any permutation $s$ acting on the set $\{a_1,\dots,a_p\}$ containing zero,
\[
{}_pF_q
\left[
\begin{array}
{r}s(a_1),\dots,s(a_p);\\
b_1,\dots,b_q;
\end{array}
z\right] = 1\,,
\]
we obtain from~\eqref{sr1} that
\begin{equation}\label{a1}
\widetilde{\ker}(u;n,0,\sigma,0) = \frac{{(-1)}^{\frac{n}2}\mathbf{i}\Gamma(n-1)}{\pi{(2u)}^{\frac{n+\sigma}2}\Gamma\left(\frac{n}2\right)}
\sum\limits_{i = 0}^\infty\frac{(-\sigma)_i\left(1-\sigma-\frac{n}2\right)_i}{i!}
W_{\frac\sigma2-i,\frac{\sigma+n-1}2}(2u)\,.
\end{equation}
Also, in view of the identity
\begin{equation}\label{impid}
\frac\pi{\sin(\pi z)} = \Gamma(z)\Gamma(1-z)\,,
\end{equation}
we have
\[
\widetilde{\ker}(u;n,0,\sigma,0) = \frac{{(-1)}^{\frac{n}2}\mathbf{i}(n-2)!}{2\pi}
\int\limits_1^{\infty}P_{\frac{n}2+\sigma-1}^{1-\frac{n}2}(t)
P_{-\frac{n}2}^{1-\frac{n}2}(t)\exp(-ut)\,\mathrm{d}t\,.
\]
From~\eqref{ALF-1st}, we get
\begin{align*}
P_{-\frac{n}2}^{1-\frac{n}2}(t) &= \Bigl\{\Gamma\Bigl(\frac{n}2\Bigr)\Bigr\}^{-1}\Bigl{(\frac{t+1}{t-1}\Bigr)}^{\frac12-\frac{n}4}
{}_1F_0\Bigl(1-\frac{n}2;-;\frac{1-t}2\Bigr)\\
&= 2^{1-\frac{n}2}\Bigl\{\Gamma\Bigl(\frac{n}2\Bigr)\Bigr\}^{-1}
{(t^2 -1)}^{\frac{n}4-\frac12}\,.
\end{align*}
Therefore,
\[
\widetilde{\ker}(u;n,0,\sigma,0) = \frac{{(-1)}^{\frac{n}2}\mathbf{i}(n-2)!}{2^{\frac{n}2}\pi\Gamma\left(\frac{n}2\right)}\int\limits_1^{\infty} P_{\frac{n}2+\sigma-1}^{1-\frac{n}2}(t)
{(t^2 -1)}^{\frac{n}4-\frac12}\exp(-ut)\,\mathrm{d}t.
\]
Using a known integral formula (see, e.g.,~\cite[Entry 2.17.7.-5]{v3})
 \begin{equation}\label{KIF-3}
  \int\limits_b^{\infty}{(t^2 -b^2)}^{-\frac{\mu}2}\exp(-pt)P_\nu^\mu\Bigl(\frac{t}b\Bigr)\,\mathrm{d}t = \Bigl{(\frac{2b}\pi\Bigr)}^{\frac12}p^{\mu-\frac12}K_{\nu+\frac12}(bp)
\end{equation}
\begin{equation*}
   (b \in \mathbb{R}^+,\Re(p)>0, \Re(\mu)<1)\,,
\end{equation*}
we get
\begin{equation}\label{a2}
\widetilde{\ker}(u;n,0,\sigma,0) = \frac{{(-1)}^{\frac{n}2}\mathbf{i}{(2u)}^{\frac{1-n}2}(n-2)!}{\pi^{\frac32}\Gamma\left(\frac{n}2\right)}K_{\frac{n-1}2+\sigma}(u)\,.
\end{equation}
Finally, equating~\eqref{a1} and~\eqref{a2} leads to the desired
identity~\eqref{thm1-sec3-eq1}.
\end{proof}

The particular case $n = 2$ of~\eqref{thm1-sec3-eq1} gives
\begin{equation*}\label{thm1-sec3-eq1:n = 2}
 K_{\sigma+\frac12}(u) = {(2u)}^{-\frac{\sigma+1}2}\pi^{\frac12}\sum\limits_{j = 0}^\infty
\frac{[(-\sigma)_j]^2}{j!}W_{\frac\sigma2-j,\frac{\sigma+1}2}(2u)\,.
\end{equation*}

\section{A series involving parabolic cylinder functions}\label{sec4}

A series associated with parabolic cylinder functions can be evaluated as in the following theorem.

\begin{theorem}\label{thm2-sec4} Let $\Re(u)>0$. Then
\begin{equation}\label{thm2-sec4-eq1}
    \sum_{j = 0}^{\infty}\,2^{2j}\, \Gamma \biggl(j+ \frac{1}{2}\biggr)
   D_{-2j-1} (2\sqrt{u}) = \frac{\sqrt{\pi}\,e^{-2u}}{2\,\sqrt{u}}\,.
\end{equation}
\end{theorem}

\begin{proof}
Taking the definition
\[
\Gamma(z) = \int\limits_0^\infty t^{z-1}\,e^{-t}\,\mathrm{d}t \quad \left(\Re(z)>0\right)
\]
of the Gamma function,
from~\eqref{ir} we obtain
\begin{equation}\label{c1}
\aligned \widetilde{\ker}\Bigl(u; 1,0,-\frac12,0\Bigr) &
= \frac1{2{\bf
i}}\int\limits_1^{\infty}\,P_{-1}^{\frac12}(t)\,P_{-\frac12}^{\frac12}(t)\,\exp(-u
t)\,\mathrm{d}t\\
&= \frac1{\sqrt{2}\pi\mathbf{i}}\int\limits_1^{\infty}{(t-1)}^{-\frac12}\,\exp(-u t)\,\mathrm{d}t = \frac1{\mathbf{i} \sqrt{2\pi u}\exp(u)}\,.
\endaligned
\end{equation}

 From the proof of Lemma~\ref{lem-Ker} and~\eqref{c} we get
\begin{equation}\label{c2}
\aligned \widetilde{\ker}\Bigl(u; 1,0,-\frac12,0\Bigr)&
= \frac{\sqrt{2}}{\pi\mathbf{i}}\sum\limits_{j = 0}^\infty2^\mathbf{i}\int\limits_{0}^{\infty}\frac{t^{j-\frac12}\,\mathrm{d}t}{{(t+2)}^{j+1}\exp(u t)}\\
& = \frac{\sqrt{2}\,e^{u}}{\pi\mathbf{i}}\,\sum_{j = 0}^{\infty}2^{2j} \Gamma \biggl(j+
\frac{1}{2}\biggr)
   D_{-2j-1} \bigl(\sqrt{4u}\bigr)\,.
\endaligned
\end{equation}
Now, equating~\eqref{c1} and~\eqref{c2} yields the desired identity~\eqref{thm2-sec4-eq1}.
\end{proof}

\section{Orthogonality relations for kernels $\widetilde{\ker}(u;n,k_0,\sigma,\hat\sigma)$}

It is noted that some particular cases of the kernel
$\widetilde{\ker}(u;n,k_0,\sigma,\hat\sigma)$ coincide with
known variants of Bessel functions. For example, applying~\eqref{impid} and
two known integral formulas (see~\cite[Entries 2.3.5.-4 and
2.3.5.-5]{v1})
\begin{equation*}\label{IF-235-4}
  \int\limits_{a}^{\infty}\, {\left(t^2 -a^2 \right)}^{\beta-1}e^{-p t}\, \mathrm{d}t
     = \frac{\Gamma (\beta)}{\sqrt{\pi}}{\left(\frac{2a}{p}\right)}^{\beta-\frac{1}{2}} K_{\beta-\frac{1}{2}} (ap)
 \end{equation*}
\begin{equation*}
  \left(\min \{\Re(\beta), \Re(p)\}>0\right)
\end{equation*}
and
\begin{equation*}\label{IF-235-5}
  \int\limits_{a}^{\infty} e^{\varepsilon \mathbf{i} \lambda t} {\left(t^2 -a^2 \right)}^{\beta-1}\, \mathrm{d}t
     = \varepsilon \mathbf{i}\sqrt{\pi}\Gamma (\beta)\,2^{\beta-\frac{3}{2}}{\left(\frac{\lambda}{a}\right)}^{\frac{1}{2}-\beta} H^{\left(\frac{3}{2}-\frac{\varepsilon}{2}\right)}_{\frac{1}{2}-\beta} (a \lambda)
 \end{equation*}
\begin{equation*}
  \left[\varepsilon = \pm 1;\Re(\beta)>0,\pm \Re(\textbf{i} \lambda)<0
  \quad \left(0<\Re(\beta)<1, \Re(\mathbf{i} \lambda) = 0\right) \right]
\end{equation*}
to~\eqref{ir}, respectively, we obtain, for even $n$,
\begin{equation*}\label{a3}
 \widetilde{\ker}(u;n,0,0,0) = {(-1)}^{\frac{n}2}{(2u)}^{\frac{1-n}2}\pi^{-\frac32}
\mathbf{i}(n-2)!\left[\Gamma\left(\frac{n}2\right)\right]^{-1}
K_{\frac{n-1}2}(u) \quad (\Re(u)>0)
\end{equation*}
and
\begin{equation*}\label{KHL-a4}
 \widetilde{\ker}(\varepsilon\mathbf{i}\omega;n,0,0,0) = {(-1)}^{\frac{n}2}2^{-\frac{n+1}2}\omega^{\frac{1-n}2}\pi^{-\frac12}\mathbf{i}^{\frac{n+1}2}(n-1)!\left[\Gamma\left(\frac{n}2\right)\right]^{-1}H_{\frac{1-n}2}^{(1+\mathbf{\eta})}(\omega)
\end{equation*}
\begin{equation*}
  \Bigl({(-1)}^\mathbf{\eta}\Re(\mathbf{i}\omega)>0\Bigr)\,.
\end{equation*}

 Therefore,
the kernel $\widetilde{\ker}(u;n,k_0,\sigma,\hat\sigma)$ as a
function can be considered as a generalization of some functions
associated with the Bessel functions. Certain properties of
$\widetilde{\ker}(u;n,k_0,\sigma,\hat\sigma)$ become
generalizations and analogues of those associated Bessel
functions. Here, in order to prove a family of orthogonality
relations for $\widetilde{\ker}(u;n,0,\sigma,0)$ and
$\widetilde{\ker}(u;n,0,0,\hat\sigma)$, we use a known
integral formula belonging to the Kontorovich-Lebedev integral
transform (see, e.g.,~\cite[Eq. (1.1)]{gm} and two references
therein include the origin of this transformation)
\begin{equation}\label{or}
\int\limits_0^{\infty}u^{-1}K_{\mathbf{i}\rho}(u)K_{\mathbf{i}\hat\rho}(u)\,\mathrm{d}u = \frac{\pi^2 \,\delta(\rho-\hat\rho)}{2\rho\sinh(\pi\rho)}
\quad (\min \left\{\Re(\rho),\Re(\hat\rho) \right\}>0),
\end{equation}
where $\delta$ is the Dirac delta function.
In~\cite{psky}, several approaches of proof of~\eqref{or} were given.
In~\cite{sb},~\eqref{or} was proved in a simpler way than those in~\cite{psky}
by appealing to a technique occasionally used in mathematical physics.

Applying the the relation between Hankel and MacDonald functions (see \cite[Eq. 5.33]{m})
\begin{equation*}\label{H-M-R}
 H_\nu^{(1+\mathbf{\eta})}(u) = {(-1)}^\mathbf{\eta}
2{(\pi\mathbf{i})}^{-1}\,
\exp\biggl(\frac{{(-1)}^{\mathbf{\eta}+1}\nu\pi\mathbf{i}}2\biggr)
K_\nu\biggl[\exp\biggl(\frac{{(-1)}^{\mathbf{\eta}+1}\pi\mathbf{i}}
2\biggr)u\biggr]
\end{equation*}
to~\eqref{or}, we get an orthogonality formula for Hankel functions
\begin{equation*}\label{Ortho-HFts}
  \int\limits_0^{\infty}u^{-1}H_{\mathbf{i}\rho}^{(1+\mathbf{\eta})}(u)H_{\mathbf{i}\hat\rho}^{(1+\mathbf{\eta})}(u)\,\mathrm{d}u = -\frac{2\exp\big({(-1)}^\mathbf{\eta}\pi\rho\big)\delta(\rho-\hat\rho)}
{\rho\sinh(\pi\rho)}
\end{equation*}
\begin{equation*}
 (\min \left\{\Re(\rho),\Re(\hat\rho) \right\}>0).
\end{equation*}
In~\cite{gm}, the Kontorovich-Lebedev transform with Hankel function as a kernel
was discussed in a detailed manner.

Assume here that $u\in\mathbb{R}$. The orthogonality formulae for
the kernel functions are given in the following theorems.

\begin{theorem}\label{thm3-sec5} The following integral formula holds.
   \begin{multline}\label{thm3-sec5-eq1}
\int\limits_0^{\infty}\widetilde{\ker}\Bigl(u;2,0,-\frac12+\mathbf{i}\rho,0\Bigr)\widetilde{\ker}\Bigl(u;2,0,-\frac12+{\bf
i}\hat\rho,0\Bigr)\,\mathrm{d}u\\
= -\frac{\delta(\rho-\hat\rho)}{4\pi\rho\,\sinh(\pi\rho)} \quad
\left(\min \left\{\Re(\rho),\,\Re(\hat\rho) \right\}>0\right).
\end{multline}
\end{theorem}

\begin{proof}
In view of~\eqref{a2}, for any even $n$ we have
\begin{multline}\label{thm3-sec5-pf1}
 \widetilde{\ker}\Bigl(u;n,0,\frac{1-n}2+\mathbf{i}\rho,0\Bigr)\widetilde{\ker}\Bigl(u;n,0,\frac{1-n}2+\mathbf{i}\hat\rho,0\Bigr)\\
 =
-\frac{[(n-2)!]^2}{{(2u)}^{n-1}\pi\bigl(\frac{n}2-1\bigr)!}
K_{\mathbf{i}\rho}(u)K_{\mathbf{i}\hat\rho}(u)\,.
\end{multline}
Then, setting here $n = 2$ and integrating both sides of~\eqref{thm3-sec5-pf1} with respect to $u$ from~$0$ to $\infty$
with the aid of~\eqref{or}
 gives~\eqref{thm3-sec5-eq1}.
\end{proof}

\begin{theorem}\label{thm4-sec5} The following integral formula holds.
   \begin{multline*}
 \int\limits_0^{\infty}\widetilde{\ker}\Bigl(u;2,0,0,-\frac12+\mathbf{i}\rho\Bigr)\widetilde{\ker}\Bigl(u;2,0,0,-\frac12+\mathbf{i}\hat\rho\Bigr)\,\mathrm{d}u\\
 = -\frac{\bigl{(-\frac12+\mathbf{i}\rho\bigr)}^2 \pi\sinh(\pi\rho)}{4\rho\cosh^2 (\pi\rho)}\delta(\rho-\hat\rho)
\quad (\min \{\Re(\rho),\Re(\hat\rho) \}>0
).
\end{multline*}
\end{theorem}

\begin{proof}
We derive from~\eqref{ALF-1st} that
$P_\nu^{-\nu}(z) = 2^{-\nu}\{\Gamma(\nu+1)\}^{-1}{(z^2 -1)}^{\frac\nu2}$.
Applying it to~\eqref{ir}, we obtain
\[
\widetilde{\ker}(u;n,0,0,\hat\sigma) = \frac{{(-1)}^{\frac{n}2}\mathbf{i}(\hat\sigma)_{n-1}}
{2^{\frac{n}2}\tan(\pi\hat\sigma)\Gamma\bigl(\frac{n}2\bigr)}
\int\limits_0^{\infty}P_{-\frac{n}2-\hat\sigma}^{1-\frac{n}2}(t){(t^2 -1)}^{\frac{n}4-\frac12}
e^{-ut}\,\mathrm{d}t\,.
\]
Using~\eqref{KIF-3}, we get
\[
\widetilde{\mathrm{ker}}(u;n,0,0,\hat\sigma) = \frac{{(-1)}^{\frac{n}2}\mathbf{i}\cot(\pi\hat\sigma)(\hat\sigma)_{n-1}}
{{(2u)}^{\frac{n-1}2}\sqrt{\pi}\Gamma\bigl(\frac{n}2\bigr)}K_{\hat\sigma+\frac{n-1}2}(u)\,,
\]
and, therefore,
\begin{multline*}
\int\limits_0^{\infty}\widetilde{\ker}\left(u;2,0,0,-\frac12+{\bf
i}\rho\right)\,\widetilde{\ker}\left(u;2,0,0,-\frac12+{\bf
i}\hat\rho\right)\,\mathrm{d}u\\
= -\frac{\bigl(-\frac12+\mathbf{i}\rho\bigr)\bigl(-\frac12+\mathbf{i}\hat\rho\bigr)}{2\pi
\coth(\pi\rho)\coth(\pi\hat\rho)}\int\limits_0^{\infty}u^{-1}K_{\mathbf{i}\rho}(u)K_{\mathbf{i}\hat\rho}(u)\,\mathrm{d}u.
\end{multline*}
Applying~\eqref{or}, we complete the proof.
\end{proof}

\begin{theorem}\label{thm5-sec5} The following integral formula holds.
\begin{multline}\label{thm5-sec5-eq1}
\int\limits_0^{\infty}\widetilde{\ker}\left(u;2,0,-\frac12+\mathbf{i}\rho,0\right)\widetilde{\ker}\Bigl(u;2,0,0,-\frac12+\mathbf{i}\hat\rho\Bigr)
\,\mathrm{d}u\\
= -\frac{\mathbf{i}\cosh(\pi\rho)\delta(\rho-\hat\rho)}{4\rho\sinh^2 (\pi\rho)}
\quad (\min\{\Re(\rho),\Re(\hat\rho) \}>0).
\end{multline}
\end{theorem}

\begin{proof}
Integrating both sides of the equality
\begin{equation*}
\widetilde{\ker}\Bigl(u;2,0,-\frac12+\mathbf{i}\rho,0\Bigr)\widetilde{\ker}\Bigl(u;2,0,0,-\frac12+\mathbf{i}\hat\rho\Bigr) =
\frac{\mathbf{i}\coth(\pi\rho)}{2\pi^2}K_{\mathbf{i}\rho}(u)K_{\mathbf{i}\hat\rho}(u)
\end{equation*}
with respect to
$u$ from $0$ to $\infty$ with the aid of~\eqref{or}, we derive~\eqref{thm5-sec5-eq1}.
\end{proof}

\begin{theorem}\label{thm6-sec5} The following integral formula holds.
\begin{multline*}
\int\limits_0^{\infty}\widetilde{\ker}(u;1,0,0,\mathbf{i}\rho)\widetilde{{\ker}}(u;3,0,0,-1+\mathbf{i}\hat\rho)\,\mathrm{d}u\\
= -\frac{\rho+\mathbf{i}}{2\sinh(\pi\rho)}\delta(\rho-\hat\rho)
\quad (\min \{\Re(\rho),\Re(\hat\rho) \}>0).
\end{multline*}
\end{theorem}

\begin{proof}
We obtain from~\eqref{ALF-1st} that for any odd $n$
\begin{equation}\label{kf1}
\aligned
\widetilde{{\ker}}(u;n,0,0,\hat\sigma) &= \frac{{(-1)}^{\frac{n-1}2}(\hat\sigma)_{n-1}}{2\mathbf{i}}\int\limits_1^{\infty}P_{\frac{n}2-1}^{1-\frac{n}2}(t)
P_{-\frac{n}2-\hat\sigma}^{1-\frac{n}2}(t)e^{-ut}\,\mathrm{d}t\\
& = \frac{{(-1)}^{\frac{n+1}2}\mathbf{i}(\hat\sigma)_{n-1}}{{(2u)}^{\frac{n-1}2}\sqrt{\pi}\Gamma\bigl(\frac{n}2\bigr)}
K_{\hat\sigma+\frac{n-1}2}(u)\,.
\endaligned
\end{equation}
Integrating both sides of
\[
\widetilde{{\ker}}(u;1,0,0,\mathbf{i}\rho)\widetilde{{\ker}}(u;3,0,0,-1+\mathbf{i}\hat\rho) = -\frac{\hat\rho^2 +\mathbf{i}\hat\rho}{\pi^2 u}
K_{\mathbf{i}\rho}(u)K_{\mathbf{i}\hat\rho}(u)
\]
and using~\eqref{or},
we complete the proof.
\end{proof}

\begin{theorem}\label{thm7-sec5} The following integral formula holds.
   \begin{multline}\label{thm7-sec5-eq1}
\int\limits_0^{\infty}\widetilde{{\ker}}\left(u;1,0,\mathbf{i}\rho,0\right)\,\widetilde{{\ker}}\left(u;3,0,0,-1+\mathbf{i}\hat\rho\right)\,\mathrm{d}u\\
= -\frac{(\rho+\mathbf{i})\delta(\rho-\hat\rho)}{2\sinh(\pi\rho)} \quad (\min
\{\Re(\rho),\Re(\hat\rho) \}>0).
\end{multline}
\end{theorem}

\begin{proof}
For any odd $n$, we have
\begin{equation}\label{kf2}
\aligned
\widetilde{{\ker}}(u;n,0,\sigma,0) &= \frac{{(-1)}^{\frac{n-1}2}(0)_{n-1}}{2\mathbf{i}}
\int\limits_1^{\infty}P_{\frac{n}2+\sigma-1}^{1-\frac{n}2}(t)
P_{-\frac{n}2}^{1-\frac{n}2}(t)e^{-ut}\,\mathrm{d}t\\
& = \frac{{(-1)}^{\frac{n+1}2}\mathbf{i}(0)_{n-1}}{{(2u)}^{\frac{n-1}2}\sqrt{\pi}\Gamma\bigl(\frac{n}2\bigr)}
K_{\sigma+\frac{n-1}2}(u)\,.
\endaligned
\end{equation}
The equality~\eqref{thm7-sec5-eq1} follows from the identity
\begin{equation*}
\widetilde{{\ker}}(u;1,0,\mathbf{i}\rho,0)\widetilde{{\ker}}(u;3,0,0,-1+\mathbf{i}\hat\rho) = \frac{\mathbf{i}\hat\rho(-1+\mathbf{i}\hat\rho)}{\pi u}K_{\mathbf{i}\rho}(u)K_{\mathbf{i}\hat\rho}(u).
\qedhere
\end{equation*}
\end{proof}

\section{Concluding Remarks}\label{sec6}

In this paper we have shown that the kernel $\widetilde{{\ker}}(u;n,k_0,\sigma,\hat\sigma)$ plays the same role as the
function $P_{mn}^l$ and the kernel $K(w,z,g)$ in~\cite{v3} in the
sense that some particular cases of $P_{mn}^l$ reduce to
 Jacobi and Legendre polynomials, Legendre and Bessel
functions, the kernel $K(w,z,g)$ can be expressed in terms of
gamma function. Here, properties of the kernel function
$\widetilde{\ker}(u;n,k_0,\sigma,\hat\sigma)$ yield those
identities
 corresponding to variants of Bessel functions and their related functions. For example, choosing
$\sigma = 0$ in Theorem~\ref{thm1-sec3}, we obtain the well-known relation~\cite[Entry~9.235.2]{gr}
\[
K_{\frac{n-1}2}(u) = \sqrt{\frac\pi{2u}}\,W_{0,\frac{n-1}2}(2u)\,.
\]
For other examples, the particular cases of of~\eqref{kf1} and~\eqref{kf2} when $\hat\rho = \rho$ give the following integral
formula:
 \[
\int\limits_0^\infty\widetilde{{\ker}}^2 (u;1,0,\mathbf{i}\rho,0)\,\mathrm{d}u = \int\limits_0^\infty\widetilde{\ker}^2 (u;1,0,0,\mathbf{i}\rho)\,\mathrm{d}u = -\frac{\mathrm{sech}(\pi\rho)}{16\pi^2}\,,
\]
which may be considered as an analogue of the following known
integral formula~\cite[Lemma~2.3]{cr}
\begin{equation*}\label{anal-sec6-eq1}
 \int\limits_0^{\infty}[K_{\mathbf{i}\rho}(2\pi
u)]^2 \,\mathrm{d}u = \frac\pi{8\cosh(\pi\rho)}.
\end{equation*}

\EditInfo{11 April 2020}{10 June 2020}{Karl Dilcher}

\end{paper}